\documentclass[]{article}
\usepackage{a4,amssymb}

\def\proofende{\hfill\rule{2.2mm}{2.2mm}\vspace{1ex}}
\def\re{{\rm Re}\,}

\DeclareMathAlphabet{\bmi}{OML}{cmm}{b}{it}

\newtheorem{satz}{Satz}[section]
\newtheorem{lemma}[satz]{Lemma}

\newtheorem{definition}[satz]{Definition}
\newtheorem{theorem}[satz]{Theorem}
\newtheorem{corollary}[satz]{Corollary}

\newtheorem{example}[satz]{Example}

\title{A Review on Realization Theory for Infinite-Dimensional
Systems\footnote{This paper was
supported by the Volkswagen Stiftung (RiP program at
Oberwolfach) and by the Deutsche Forschungsgemeinschaft.}}

\author{Birgit Jacob\\ 
Bergische Universit\"{a}t Wuppertal\\
Fakult\"{a}t f\"{u}r Mathematik und Naturwissenschaften\\
Arbeitsgruppe Funktionalanalysis\\
Gau\ss stra\ss e 20,
D-42119 Wuppertal, Germany\\
jacob@math.uni-wuppertal.de
\and
Hans Zwart\\ Department of Applied Mathematics\\ University of Twente\\
P.O.~Box 217, 7500 AE Enschede\\
The Netherlands\\
h.j.zwart@utwente.nl}

\begin{document}
\maketitle

\begin{abstract}
We give an introduction to the realization theory for infinite-dimensional systems. That is, we show that for any function $G$, analytic and bounded in the right half of the complex plane, there exist (unbounded) operators $A,B,C$ such that $G(s_1)-G(s_2) = (s_2-s_1) C(s_1 I-A)^{-1}(s_2 I-A)^{-1}B$. Here $A$ is the infinitesimal generator of a strongly continuous semigroup on a Hilbert space, and $B$ and $C$ are admissible input and output operators, respectively. Our results summarise and clarify the results as found in the literature, starting more than 40 years ago.
\end{abstract}

\section{Introduction}

Already since the beginning of infinite-dimensional systems theory,
there has been interest in the state-space realization problem.
The state-space realization problem is the problem of finding, for
a given function, a system in state-space form whose transfer function
equals this given function.
If we start with a rational function, then it is well-known that
we can always find a system with a
finite-dimensional state-space whose transfer function is the given
rational function. Moreover, it is always possible to find a controllable
and observable realization and all controllable and observable
realizations are equivalent, i.e., let $(A_1,B_1,C_1,D_1)$ and
$(A_2,B_2,C_2,D_2)$ be two realizations  of the transfer function $G(s)$ that are controllable and
observable, then there exists an invertible matrix $S$ such that
$A_1= S A_2 S^{-1}, B_1= S B_2, C_1 = C_2 S^{-1}$, and $D_1=D_2$.

Since every finite-dimensional system has a
rational transfer function, for a non-rational transfer function it is only possible to find
a (state-space) realization with an infinite-dimensional state space.
For functions that are analytic and bounded in some right-half
plane, the realization problem was investigated by a number of
people, e.g.\ Baras and Brockett \cite{babr73}, Fuhrmann \cite{fu81},
Helton \cite{he76}, Yamamoto \cite{ya81,ya82},
Salamon \cite{sa89} and  Weiss \cite{we89b,we97}.
Here we present the realization theory in the language of well-posed
linear systems.

We end this section with some notation and well-known results, see e.g.\ \cite[Appendix A.6]{CuZw95}, Duren \cite{Dure00} or Rudin \cite{Rudi87}.
\begin{eqnarray*}
\mathbb C_\delta^+ &:=& \{ z\in\mathbb C \mid \re z >\delta\},\quad
\delta\in \mathbb R,\\
\mathbb C_\delta^- &:=& \{ z\in\mathbb C \mid \re z <\delta\},\quad
\delta\in \mathbb R,\\
\bmi{H}_\infty(\Omega)  &:=& \{f:\Omega\rightarrow \mathbb C \mid f
\mbox{ is holomorphic and bounded}\},\quad \Omega\subset \mathbb C,\\
\bmi{H}_2({\mathbb C}_0^+)  &:=& \{f:{\mathbb C}_0^+ \rightarrow \mathbb C \mid f
\mbox{ holomorphic and }\sup_{r >0} \int_{-\infty}^{\infty} |f(r+i \omega)|^2 d \omega < \infty \},\\
\bmi{H}_2({\mathbb C}_0^-)  &:=& \{f:{\mathbb C}_0^- \rightarrow \mathbb C \mid f
\mbox{ holomorphic and }\sup_{r < 0} \int_{-\infty}^{\infty} |f(r+i \omega)|^2 d \omega < \infty \}.
\end{eqnarray*}

Furthermore, a holomorphic function $G:\mathbb C_0^+ \rightarrow
\mathbb C$ is called {\em inner} if $|G(z)|\le 1$ for $z\in
\mathbb C_0^+$, and $|G(it)|=1$ for almost every $t\in \mathbb R$.

Clearly, every inner function is an element of
$\bmi{H}_\infty(\mathbb C_0^+)$. Hence we may define
$G(i\cdot)$ as the non-tangential limits of $G(z)$. This limit
exists almost everywhere, and so the condition, $|G(it)|=1$ a.e., makes sense.

On $\bmi{H}_2({\mathbb C}_0^+)$ and $\bmi{H}_2({\mathbb C}_0^-)$ we define the following inner product from $\bmi{L}_2(i{\mathbb R})$;
\begin{equation}
\label{eq:1}
  \langle f, g \rangle = \frac{1}{2\pi} \int_{-\infty}^{\infty} f(i \omega) \overline{g(i\omega)} d \omega.
\end{equation}
Then the following holds
\begin{enumerate}
\item
  $\bmi{H}_2({\mathbb C}_0^+)^{\perp}= \bmi{H}_2({\mathbb C}_0^-)$ ;
\item
  $\bmi{H}_2({\mathbb C}_0^+) \oplus \bmi{H}_2({\mathbb C}_0^-) = \bmi{L}_2(i{\mathbb R})$.
\item
  The Laplace transform is an isometric isometry between $\bmi{L}_2(0,\infty)$ and $\bmi{H}_2({\mathbb C}_0^+)$. Similarly, $\bmi{L}_2(-\infty,0)$ and $\bmi{H}_2({\mathbb C}_0^-)$ are isometrically isomorph.
\item
  The Fourier transform is an isometric isomorphism between $\bmi{L}_2({\mathbb R})$ and $\bmi{L}_2(i{\mathbb R})$.
\end{enumerate}
The latter two results are known as Paley-Wiener theorem. 

\section{Well-posed linear systems and realization theory}
\label{reali}

A quite large and well-studied class of infinite-dimensional linear systems
is the class of well-posed linear systems introduced by Salamon
\cite{sa89} and Weiss \cite{we89}, \cite{we89c}. By now there are excellent books on well-posed linear systems, see e.g.\ \cite{Staf05} or \cite{TuWe09}.

Let $T(t)$ be a $C_0$-semigroup on the separable Hilbert space $H$, and
let $A$ its generator.
We define the
space $H_{-1}$ to be the completion of $H$ with respect to the norm
\[ \|x\|_{-1} := \|(\beta I- A)^{-1} x\|\]
and the space $H_1$ to be $D(A)$ with the norm
\[ \|x\|_1 := \|(\beta I- A) x\|,\]
where $\beta\in\rho(A)$, the resolvent set of $A$.
It is easy to verify that the topology of $H_{-1}$ and $H_1$ does not depend on
$\beta\in\rho(A)$. Moreover, $\|\cdot\|_1$ is equivalent to the graph
norm on $D(A)$, so $H_1$ is complete. In Weiss \cite[Remark 3.4]{we89} it
is shown that $T(t)$ has a restriction to a $C_0$-semigroup on $H_1$ whose
generator is the restriction of $A$ to $D(A)$, and $T(t)$ can be extended
to a $C_0$-semigroup on $H_{-1}$ whose generator is an extension of $A$ with
domain $H$. Therefore, we get
\[ A\in {\cal L}(H_1,H)\quad\mbox{ and }\quad  A\in {\cal
  L}(H,H_{-1}).\]
$H_{-1}$ equals the dual of $D(A^\ast)$, where we have equipped
$D(A^\ast)$ with the graph norm (see \cite{we89}).
Following \cite{we89} and \cite{we89c} we introduce admissible control operators
and observation operators for $T(t)$.
\begin{definition}
\begin{enumerate}
\item Let $B\in {\cal L}(\mathbb{C},H_{-1})=H_{-1}$. For $t\ge 0$ we
define the operator ${\cal B}_t: \bmi{L}_2(0,\infty)\rightarrow
H_{-1}$ by
\[ {\cal B}_t u:= \int_0^t T(t-\rho)B u(\rho)\,d\rho.\]
Then $B$ is called an {\em admissible control operator for
$T(t)$}, if for some (and hence any)
$t> 0$, ${\cal B}_t\in{\cal L}(\bmi{L}_2(0,\infty), H)$.
\item Let $B$ be an admissible control operator for $T(t)$.
  $B$ is called an {\em infinite-time admissible control operator for
  $T(t)$}, if $T(\cdot)Bu(\cdot):[0,\infty)\rightarrow H_{-1}$ is integrable for every
  $u\in \bmi{L}_2(0,\infty)$, and the operator ${\cal B}_{\infty}:\bmi{L}_2(0,\infty)
  \rightarrow H_{-1}$, given by
  \[ {\cal B}_{\infty} u := \int_0^{\infty} T(t) B u(t) dt, \]
  satisfies ${\cal B}_{\infty} \in {\cal L}(\bmi{L}_2(0,\infty),H)$.
\item Let $C\in {\cal L}(H_1,\mathbb{C})$. Then $C$ is
called an {\em admissible observation
operator for $T(t)$}, if for
some (and hence any) $t> 0$, there is some $K>0$ such that
\[ \|CT(\cdot)x\|_{\bmi{L}_2(0,t)}\le K\|x\|, \quad x\in D(A).\]
\item Let $C$ be an admissible observation operator for $T(t)$.
We call $C$ an {\em infinite-time admissible observation operator}\/
if there is some $K>0$ such that
\[ \|CT(\cdot)x\|_{\bmi{L}_2(0,\infty)}\le K\|x\|, \quad x\in D(A).\]
\end{enumerate}
\end{definition}

By definition, every infinite-time admissible control operator for $T(t)$ is
an admissible control operator for $T(t)$. If $T(t)$ is exponentially stable,
then the two notions of admissibility coincide. Similar statements hold for admissible
observation operators for $T(t)$. Moreover, $B$ is an (infinite-time) admissible
control operator for $T(t)$ if and only if $B^*$ is an (infinite-time) admissible
observation operator for $T^*(t)$.

Let $C$ be an admissible observation operator for $T(t)$. Then for
$t\ge 0$ the operator $\Psi_t\in{\cal L}(D(A), \bmi{L}_2(0,t))$,
given by $\Psi_t x:= CT(\cdot)x$, has a unique extension to
${\cal L}(H, \bmi{L}_2(0,t))$ (again denoted by $\Psi_t$).
Similarly, if $C$ is an infinite-time admissible observation, then
we can extend the operator $\Psi_\infty\in {\cal L}(D(A), \bmi{L}_2(0
,\infty))$, given by $\Psi_\infty x:= CT(\cdot)x$, to $ {\cal L}(H, \bmi{L}_2(0
,\infty))$. 
\begin{definition}\label{def1i}
Let $B$ be an admissible control operator for $T(t)$, then
\begin{enumerate}
\item
  {\em $(T,B)$ is exactly controllable in finite time} if there exists
  a time $t_0$ such that
  $ {\rm Im}\, {\cal B}_{t_0} = H$.
\item
  {\em $(T,B)$ is approximately controllable} if
  $\overline{\cup_{t\geq 0} {\rm Im}\,{\cal B}_t} = H$.
\item
  {\em $(T,B)$ is exactly controllable} if $B$ be an
  infinite-time admissible control operator for $T(t)$ and
   ${\rm Im}\, {\cal B}_{\infty} = H$.
\end{enumerate}
Let $C$ be an admissible observation operator for $T(t)$, then
\begin{enumerate}
\item
  {\em $(T,C)$ is exactly observable in finite time} if  there exists
  a time $t_0$ such that
  $\|\Psi_{t_0} x \| \geq c\|x\|$
  for some positive $c$ and every $x\in X$.
\item
  {\em $(T,C)$ is approximately observable} if
  $\cap_{t\geq 0} \ker\Psi_t = \{0\}$.
\item
  {\em $(T,C)$ is exactly observable} if  $\Psi_{\infty}
  \in {\cal L}(H,\bmi{L}_2(0,\infty))$ and
  $\|\Psi_{\infty} x \| \geq c\|x\|$
  for some positive $c$ and every $x\in X$.
\end{enumerate}
\end{definition}

If $B$ is infinite-time admissible, then it is easy to see that
$(T,B)$ is approximately controllable if and only if
$\overline{\mbox{\rm Im\,}{\cal
B}_{\infty}} = H$. A similar statement holds for the observation
operator and approximate observability.
It is easy to see, that controllability and observability are dual
notions, i.e., $(T,B)$ is approximately (exactly) controllable if and
only if $(T^*,B^*)$ is approximately (exactly) observable.
We are now in the position to introduce well-posed linear systems.

\begin{definition}\label{def1}
$(T,B,C,G)$ is called a {\em well-posed linear system}
if the following holds
\begin{enumerate}
\item $T(t)$ is a $C_0$-semigroup,
\item $B$ is an admissible control operator for $T(t)$,
\item $C$ is an admissible observation operator for $T(t)$,
\item There exists a constant $\rho>0$ such that $G\in \bmi{H}_\infty(\mathbb C_\rho^+)$ and
      \begin{equation}\label{eqn14}
        \frac{G(s)-G(z)}{z-s}= C(sI-A)^{-1}(zI-A)^{-1}B,\qquad s,z\in \mathbb C_\rho^+, s\not=z,
      \end{equation}
      where $A$ is the generator of $T(t)$.
\end{enumerate}
\end{definition}

Here $G$ is called the {\em transfer function}.
A transfer function is determined by $(T,B,C)$ up to an additive constant
operator. Note, that this definition of a well-posed linear system is not the standard one
introduced by Weiss \cite{we89} \cite{we89c}, but it is equivalent
to Weiss's definition, see Curtain and Weiss \cite{cuwe89}.
Moreover, this definition of a well-posed linear system is equivalent to Salamon's definition of
a time-invariant, linear control system, see Salamon \cite{sa89} and Weiss \cite{we94c}.
Following, Salamon \cite{sa89}, the system trajectory of a well-posed linear system $(T,B,C,G)$
is given by
\begin{eqnarray*}
x(t,x_0,u) &:=& T(t)x_0 +\int_0^t T(t-\rho) Bu(\rho)\,d\rho,\qquad t\ge0,\\
y(t,x_0,u) &:=& C(x(t,x_0,u)-(\mu I-A)^{-1}Bu(t))+G(\mu)u(t), \qquad t\ge0.
\end{eqnarray*}
Here $u\in \bmi{L}_2(0,\infty)$ denotes the input of the system, $x_0\in H$ denotes the
initial state, $x(t,x_0,u)$ denotes the state of the system at time $t$, and $y(t,x_0,u)$ denotes
the output at time $t$. Note, that the definition of $y(t,x_0,u)$ does not depend on $\mu\in \rho(A)$,
where $\rho(A)$ denotes the resolvent set of $A$. Moreover, if additionally $(T,B,C,G)$
is a {\em regular system}, i.e., $D:=\lim_{s\rightarrow+\infty, s \in
{\mathbb R}}G(s)$ exists, then we are able
to choose $D$ as the feedthrough term and $y(t,x_0,u)$ is given by
\[ y(t,x_0,u) = Cx(t,x_0,u)+ Du(t), \qquad t\ge0.\]
Note, that in general  $D$ does not exists.
\begin{definition}
A well-posed linear system $(T,B,C,G)$ is called
{\em exponentially stable}, if $T$ is an exponentially
stable $C_0$-semigroup. The system is called {\em infinite-time
admissible} if $B$ is an infinite-time admissible control operator
and $C$ is an infinite-time admissible observation operator for
$T(t)$.

Furthermore, it will be called {\em approximately controllable},
{\em exactly controllable}, or {\em exactly controllable in finite-time}
if $(T,B)$ has that property. Similar properties are defined concerning observability.
\end{definition}

\begin{definition}
Let $G\in \bmi{H}_\infty(\mathbb C_\delta^+)$ for some $\delta>0$. We say that
$G$ has a realization as a well-posed linear system if there exists
a well-posed linear system $(T,B,C,G)$.
\end{definition}

Salamon \cite{sa89} proved that every $G\in \bmi{H}_\infty(\mathbb C_\delta^+)$ has a
realization as a well-posed linear system. His realization is the
well-known shift realization. The shift realization has already been
studied by many mathematicians, see for example Baras and Brockett
\cite{babr73}, Helton \cite{he76}, Fuhrmann \cite{fu81}, Yamamoto
\cite{ya81,ya82}, and Weiss \cite{we97}.
Before we state the proof, we would like to present the motivation
behind this proof.

Assume that $G$ is the Laplace transform of a continuous function $h
\in {\bmi L}_1(0,\infty) \linebreak[0]\cap \bmi{L}_2(0,\infty)$. Then finding a
state-space realization just means finding a triple $(T(t),B,C)$
such that
\begin{equation}
\label{H15}
  h(t) = CT(t)B, \qquad t \geq 0.
\end{equation}
To establish this equality the idea is very simple. We choose the 
state space $\bmi{L}_2(0,\infty)$, and define for $z \in
\bmi{L}_2(0,\infty)$ the left-shift semigroup
\begin{equation}
\label{H16}
  \left[T(t) z\right](x) = z(t+x),\quad x \geq 0.
\end{equation}
Furthermore, for a continuous functions in $z\in \bmi{L}_2(0,\infty)$
we define
\begin{equation}
\label{H17}
  C z = z(0).
\end{equation}
If we now choose $B = h$, then by combining (\ref{H16}) and
(\ref{H17}) we see that
\[
  CT(t) B = \left[T(t)h\right](0)= h(t), \quad t\geq 0.
\]
This is precisely equality (\ref{H15}), and hence we have constructed a
realization.

For transfer functions that do not have a smooth inverse Laplace
transform, there are difficulties in defining
$CT(t)B$. However, as we shall show, the realization is still
possible with these choices of $T(t), B$ and $C$. 
In order to prove
the realization, it is easier to work with the Laplace transforms of
the above object. Thus $\bmi{L}_2(0,\infty)$ becomes
$\bmi{H}_2({\mathbb C}_0^+)$, and  $T(t),B$ and $C$ become their
equivalent counterpart on this space.
\begin{theorem}\label{theoreal0}
Every function $G \in \bmi{H}_\infty(\mathbb C_0^+)$ has a realization.
\end{theorem}
{\bf Proof:}
As state space $H$ we choose
\[
  H:=\bmi{H}_2({\mathbb C}_0^+).
\]
Our $C_0$-semigroup is given by
\[
  T(t)x = P_{\bmi{H}_2(\mathbb C_0^+)}[e^{t\cdot}x(\cdot)],
  \qquad x\in H,
\]
where $P_{\bmi{H}_2(\mathbb C_0^+)}$ is the orthogonal projection
from $\bmi{L}_2(i\mathbb R)$ to $\bmi{H}_2(\mathbb C_0^+)$.
This semigroup has as infinitesimal generator
\[
  (Ax)(s)= s x(s)-\check{x}(0),\qquad x\in D(A),
\]
with
\[
  D(A) =\{ x\in H\mid s\mapsto s x(s)-\check{x}(0)\in
  \bmi{H}_2(\mathbb C_0^+)\},
\]
where $\check{x}$ denotes the inverse Laplace transform of $x$.

Furthermore, we define
\[
  B=G,
\]
and
\[
  Cx=\check{x}(0) \qquad \mbox{for } x\in D(A).
\]

Having made these choices, we still have to prove that they form a
well-posed linear system with transfer function $G$.
We begin by showing that $T(t)$ is a $C_0$-semigroup on $H$.
\begin{enumerate}
\item
  Let $S(t)$ be the right-shift semigroup on $H$, i.e.~
  \[
    S(t)x:=e^{-t\cdot} x(\cdot),\quad x\in H.
  \]
  It is easy to see that this is a $C_0$ semigroup, and that the
  adjoint of this $C_0$-semigroup equals $T(t)$. As $H$ is a  Hilbert space,  the adjoint of a
  $C_0$-semigroup is again a $C_0$-semigroup. Hence, $T(t)$ is a $C_0$-semigroup.
\item
  Now we shall derive a simple formula for the resolvent operator
  of $A$. Since the right-shift semigroup has growth bound zero, so
  has its adjoint, $T(t)$, and thus we have that the open right-half
  is contained in the resolvent set of $A$. Take $\beta \in
  {\mathbb C}_0^+$, and $x,y \in H$, then
  \begin{eqnarray*}
   \langle y, (\beta I - A)^{-1} x \rangle
   &=&
   \int_0^{\infty} \langle y, e^{-\beta t} T(t) x \rangle dt \\
   &=&
   \int_0^{\infty} e^{-\overline{\beta} t}\langle T(t)^*y, x \rangle dt \\
   &=&
   \int_0^{\infty} e^{-\overline{\beta} t}
   \frac{1}{2\pi} \int_{-\infty}^{\infty} e^{-i\omega t} y(i\omega) \overline{x(i
   \omega)} d\omega dt \\
   &=&
   \frac{1}{2\pi} \int_{-\infty}^{\infty} \int_0^{\infty}
   e^{-(\overline{\beta} + i \omega) t} y(i\omega) \overline{x(i
   \omega)} dt d\omega \\
   &=&
   \frac{1}{2\pi} \int_{-\infty}^{\infty}
   \frac{1}{\overline{\beta} + i \omega} y(i\omega) \overline{x(i
   \omega)} d\omega \\
   &=&
   \frac{1}{2\pi} \int_{-\infty}^{\infty} y(i\omega)
   \overline{ \left(\frac{x(i\omega)}{\beta - i \omega}  \right) }
   d\omega \\
   &=&
   \frac{1}{2\pi} \int_{-\infty}^{\infty} y(i\omega)
   \overline{\left(\frac{x(i\omega) - x(\beta)}{\beta - i \omega} \right) }
   d\omega \\
   &=&
   \langle y, \frac{x(\cdot)-x(\beta)}{\beta - \cdot} \rangle .
  \end{eqnarray*}
  Here we have used that $x(\beta)/(\beta - \cdot)$ is in
  $\bmi{H}_2({\mathbb C}_0^-)=\left[\bmi{H}_2({\mathbb C}_0^+) \right]^{\perp}$.
  Furthermore, is it easy to see that
  $\frac{x(\cdot)-x(\beta)}{\beta-\cdot}$ is an element of
  $\bmi{H}_2({\mathbb C}_0^+)$.
  So we conclude that
  \begin{equation}
  \label{H1}
    \left[(\beta I - A)^{-1} x\right](s) =
    \frac{x(s)-x(\beta)}{\beta - s}
    \qquad \mbox{for } s \neq \beta.
  \end{equation}
\item
  The inverse Laplace transform of the function given in
  (\ref{H1}) is given by
  \[
    \check{z}(t)=\int_t^{\infty} \check{x}(\tau) e^{-\beta \tau}
    d\tau.
  \]
  Therefore, $\check{z}$ is a continuous function on
  $[0,\infty)$, and the value in zero equals $x(\beta)$.

  From (\ref{H1}) an easy calculation shows that
  \begin{equation}\label{H2}
    (Ax)(s)= s x(s)-\check{x}(0),\qquad x\in D(A),
  \end{equation}
  with
  \begin{equation}\label{H3}
    D(A) =\{ x\in H\mid s\mapsto s x(s)-\check{x}(0)\in
    \bmi{H}_2(\mathbb C_0^+)\}.
  \end{equation}
\item
  Using (\ref{H1}) we can identify $H_{-1}=D(A^*)'$
  with the space
  \begin{equation}
  \label{H4}
    \left\{f:\mathbb C_0^+\rightarrow \mathbb C\mid
    s\mapsto \frac{f(s) -f(\beta)}{\beta-s}\in H
    \mbox{ for some }\beta\in\mathbb C_0^+ \right\}.
  \end{equation}
  Moreover, the sesquilinear form
  $\langle\cdot,\cdot\rangle_{D(A^*)\times D(A^*)'}$ is given by
  \[
    \langle f,g\rangle_{D(A^*)\times D(A^*)'}
    = \int_{-\infty}^\infty f(i\omega)\overline{g(i\omega)}\,d\omega,
    \qquad f\in D(A^*), g\in D(A^*)'.
  \]
\item
  From (\ref{H4}) it follows directly that
  $\bmi{H}_{\infty}({\mathbb C}_0^+)$ can be seen as a
  subspace of $H_{-1}$, and so $B=G$ is an element of $H_{-1}$.

  Now we show that $B$ is admissible.
  Take $y \in D(A^*)$
  \begin{eqnarray*}
   \lefteqn{ \langle y,
    \int_0^{\infty}\! T(t) B u(t) dt \rangle_{D(A^*)\times D(A^*)'} }\qquad\\
    &=&
    \int_0^{\infty} \langle y,
    T(t) B u(t) \rangle_{D(A^*)\times D(A^*)'} dt \\
    &=&
    \int_0^{\infty} \langle T(t)^* y,
    B u(t) \rangle_{D(A^*)\times D(A^*)'} dt \\
    &=&
    \int_0^{\infty} \! \frac{1}{2\pi} \int_{-\infty}^{\infty} e^{-i\omega t}
    y(i\omega) \overline{G(i\omega) u(t)} d\omega dt \\
    &=&
    \frac{1}{2\pi} \int_{-\infty}^{\infty} y(i\omega) \overline{G(i\omega)}
    \overline{\int_0^{\infty} e^{i\omega t}u(t) dt} d\omega \\
    &=&
    \frac{1}{2\pi} \int_{-\infty}^{\infty} y(i\omega) \overline{G(i\omega)}
    \overline{\hat{u}(-j\omega)} d\omega \\
    &=&
    \langle y, P_{\bmi{H}_2} \left(G(\cdot) \hat{u}(-\cdot) \right)
    \rangle.
  \end{eqnarray*}
  Since for every $u \in \bmi{L}_2(0,\infty)$, we have that
  $\hat{u}(-\cdot) \in \bmi{L}_2(i{\mathbb R})$, and since $G$ is bounded
  on the imaginary axis, we have that $G(\cdot)\hat{u}(-\cdot) \in
  \bmi{L}_2(i{\mathbb R})$. Thus we have that
  \begin{equation}
  \label{H5}
    \int_0^{\infty} T(t) B u(t) dt =
    P_{\bmi{H}_2(\mathbb C_0^+)} \left(G(\cdot) \hat{u}(-\cdot) \right)
  \end{equation}
  is well-defined for every $u \in \bmi{L}_2(0,\infty)$ with values
  in $H$. Thus $B$ is an admissible control operator for $T$.
\item
  Part 2 implies that for $x \in D(A)$, $\check{x}(0)$ is
  well-defined. This immediately proves that $C$ is a well-defined
  operator on $D(A)$. From part 3 and equation (\ref{H1}), we see that
  \begin{equation}
  \label{H6}
    C(\beta I - A)^{-1} x = x(\beta).
  \end{equation}
  Since $C(\cdot I - A)^{-1}x$ is the Laplace transform of
  $CT(\cdot)x$, we get from the equation above that
  \begin{equation}
  \label{H7}
    CT(t) x = \check{x}(t),
  \end {equation}
  for every $x \in H =\bmi{H}_2(\mathbb C_0^+)$. From Paley-Wiener
  theorem we get that the $\bmi{L}_2(0,\infty)$-norm of $\check{x}$ equals
  the $H$-norm of $x$. In other words, for every $x \in H$ we have
  that $\Psi_{\infty} x := C T(t) x \in \bmi{L}_2(0,\infty)$. Hence $C$
  is an admissible observation operator for $T$.

  Note that we even have
  \begin{equation}
  \label{H8}
    \|\Psi_{\infty} x \| = \|x\|.
  \end{equation}
  Thus $(T,C)$ is exactly observable.
\item
 We now show that $G$ is a transfer function of $(T,B,C)$.
 We have to show that (\ref{eqn14}) holds.
 Combining equation (\ref{H6}) with (\ref{H1}) gives that
 \begin{eqnarray*}
   C(sI-A)^{-1}(\beta I-A)^{-1}B
   &=&
   [(\beta I-A)^{-1}B](s)=[(\beta I-A)^{-1} G](s)\\
   &=&
   \frac{G(s)-G(\beta)}{\beta-s}.
 \end{eqnarray*}
\end{enumerate}
Thus we have constructed a realization of the
transfer function $G$.
\hfill$\blacksquare$
\bigskip

For a $G \in \bmi{H}_{\infty}({\mathbb C}_0^+)$
the Hankel operator with symbol $G$ is defined as the operator $H_G
: \bmi{L}_{2}(0,\infty) \mapsto \bmi{L}_{2}(0,\infty)$
given by
\begin{equation}
\label{H10}
  \widehat{H_G u} := P_{\bmi{H}_2(\mathbb C_0^+)}
  (G(\cdot)\hat{u}(-\cdot)),\qquad u\in \bmi{L}_2(0,\infty),
\end{equation}
where $~\hat{~}$ denotes the Laplace transform.
If $(T,B,C,G)$ is a realization of $G$, $B$ is an infinite-time
admissible control operator, and $C$ is an infinite-time
admissible observation operator for $T(t)$ we get
\[
  H_G = \Psi_{\infty} {\cal B}_{\infty}.
\]
 From the Hankel operator we can derive special results.
\begin{lemma}
\label{lemH1}
If the Hankel operator with symbol $G$
has closed range, then all infinite-time admissible, approximately
controllable and approximately observable realization  are
equivalent, i.e., if $(T_1,B_1,C_1,G)$  with state space $H_1$,
and $(T_2,B_2,C_2,G)$  with state space $H_2$ are
both infinite-time admissible, approximately controllable and
approximately observable realizations, then
there exists a bounded, invertible operator $S \in {\cal
L}(H_1,H_2)$ such that
\[
  T_2 = S T_1 S^{-1}, \quad {\cal B}_{2,\infty}= S{\cal
  B}_{1,\infty}, \quad \Psi_{2,\infty} = \Psi_{1,\infty} S^{-1}.
\]
Furthermore, if a realization of $G$ is
exactly controllable and exactly observable, then
$H_G$ has closed range.
\end{lemma}
{\bf Proof} The first part is shown in Proposition 6.2 of Ober and
Wu \cite{obwu96}.
\smallskip

Let us now assume that $G$ has an exactly controllable
and exactly observable realization. This implies, that $G$ can be
written in the form
\[  H_G=\Psi_\infty {\cal B}_\infty,
\]
where ${\cal B}_\infty\in {\cal L}(\bmi{L}_2(0,\infty),H)$, $\Psi_\infty
\in {\cal L}(H,\bmi{L}_2(0,\infty))$, ${\cal B}_\infty$ is surjective and
$\|\Psi_\infty x\| \ge c\|x\|$, for some positive $c$.
The surjectivity of ${\cal B}_\infty$ implies that the range of
$H_G$ equals the range of $\Psi_\infty$, and the open mapping theorem
shows that the range of $\Psi_\infty$ is closed. Thus $H_G$ has closed
range.
\proofende

Let us remark that the above result is not true if ${\cal
B}_{\infty}$ and/or $\Psi_{\infty}$ are not bounded operators on
$\bmi{L}_2(0,\infty)$ and $H$, respectively. An example can be found at
the end of this section.
\begin{corollary}
\label{C3.7}
Suppose $H_G$ has closed range and there exists an exactly controllable
and exactly observable realization. Then every infinite-time
admissible, approximately controllable and approximately observable
realization is exactly controllable and exactly observable.
\end{corollary}

If the transfer function is inner, then there
exists a realization which is exactly controllable and
exactly observable.
\begin{theorem}\label{theoreal}
Let $G\in \bmi{H}_\infty(\mathbb C_0^+)$  be an inner function.
Then there exists an
exactly controllable and exactly observable well-posed
linear system  $(T,B,C,G)$.

Furthermore, $G$ has an
exactly controllable and exactly observable realization with the
$C_0$-semigroup having the additional property that
\begin{enumerate}
\item it is exponentially stable if and only if $\inf_{\re z\in
      (0,\alpha)} |G(z)|>0$ for some $\alpha>0$.
\item it is a group if and only if $\inf_{\re z >
      \rho} |G(z)| >0$ for some $\rho>0$.
\end{enumerate}
\end{theorem}
{\bf Proof\/} Theorem \ref{theoreal0} shows there exists
a realization of  $G$. We use the notation as introduced in
the proof of Theorem \ref{theoreal0}.
Let $V$ denote the closed subspace of $\bmi{H}_2({\mathbb C}_0^+)$
defined as
\[
  V=[G \bmi{H}_2(\mathbb C_0^+)]^{\perp},
\]
where the orthogonal complement is taken in
$\bmi{H}_2(\mathbb C_0^+)$.
We shall show that the realization $(T|_{V},G,C|_{V})$ has the
desired properties.

We begin by showing that $V$ is $T(t)$-invariant.
\begin{enumerate}
\item
  Take an arbitrary $x \in H$ and $v \in V$, then
  \begin{eqnarray*}
    \langle G x, T(t) v \rangle 
    &=&
    \langle T(t)^* G x, v \rangle 
    =
    \langle e^{-t\cdot}G(\cdot)x(\cdot), v(\cdot) \rangle\\
    &=&
    \langle G(\cdot)e^{-t\cdot}x(\cdot), v(\cdot) \rangle
    =0,
  \end{eqnarray*}
  since $e^{-t\cdot} x \in H$, and so
  $G(\cdot)e^{-t\cdot}x(\cdot) \in V^{\perp}$.

  The set of all $w$ that can be written as $G x$ is dense in
  $V^{\perp} = \overline{G \bmi{H}_2({\mathbb C}_0^+)}$,
  thus we have that $T(t)V \subset V$.

  From this we see that $T_V(t)$ defined as the restriction to $V$
  of $T(t)$ is a $C_0$-semigroup on $V$.
\item
  Now we show that $B:=G$ is an admissible control operator
  for $T_V(t)$.

  From (\ref{H5}) we see that
  \[
    {\cal B}_{\infty}u= \int_0^{\infty} T(t) B u(t)dt =
    P_{\bmi{H}_2(\mathbb C_0^+)}
    \left(G(\cdot) \hat{u}(-\cdot) \right).
  \]
  If we can show that this expression maps into $V$, then we are
  done.
  Take an $x \in H$, and consider
  \begin{eqnarray*}
    \langle G x,{\cal B}_{\infty}u \rangle
    &=&
    \frac{1}{2\pi} \int_{-\infty}^{\infty} G(i\omega) x(i\omega)
    \overline{G(i\omega) \hat{u}(-i\omega)} d\omega\\
    &=&
    \frac{1}{2\pi} \int_{-\infty}^{\infty} x(i\omega)
    \overline{\hat{u}(-i\omega)} d\omega
    =0,
  \end{eqnarray*}
  where we have used $x \in H=\bmi{H}_2(\mathbb C_0^+)$, and
  $\hat{u}(-\cdot) \in \bmi{H}_2(\mathbb C_0^+)^{\perp}$.

  So we have shown that ${\cal B}_{\infty}u$ is orthogonal to any
  $G x$ with $x \in \bmi{H}_2(\mathbb C_0^+)$. This proves
  that ${\cal B}_{\infty}$ maps into $V$.
\item
  Since $C$ is admissible for $T(t)$ it is directly clear that $C|_V$
  defined as
  \begin{equation}
  \label{H9}
    C_V x:=\check{x}(0) = C x,\quad x \in V
  \end{equation}
  is admissible for $T_V(t)$
\item
  Combining the above results we see that we found a second realization of
  $G$. We shall now prove that it is exactly
  controllable. Note that from equations (\ref{H9}) and (\ref{H8})
  it follows that $(T_V,C_V)$ is exactly observable.
\item
  To show that the range of ${\mathcal B}_{\infty}$ is closed we prove that it is a partial isometry. That is for every $u  \in \bmi{L}_2(0,\infty)$ with $u\perp \ker  {\mathcal B}_{\infty}$ there holds
  \begin{equation}
  \label{eq:15}
    \| {\mathcal B}_{\infty} u \| = \|u\|.
  \end{equation}
  Let $\hat{v} \in  \bmi{H}_2({\mathbb C}_0^-)$, and define for $s \in {\mathbb C}_0^-$, $G^{\dagger}(s) = \overline{G(-\overline{s})}$. It is easy to see that $G^{\dagger} \in \bmi{H}_{\infty}({\mathbb C}_0^-)$ and also that $G^{\dagger} \hat{v} \in  \bmi{H}_2({\mathbb C}_0^-)$.  For $s \in {\mathbb C}_0^+$ define $\hat{q}(s) = G^{\dagger}(-s) \hat{v}(-s)$. Then $\hat{q} \in \bmi{H}_{2}({\mathbb C}_0^+)$. Finally, we denote by $q$ the inverse Laplace transform of $\hat{q}$. We claim that $q \in \ker  {\mathcal B}_{\infty}$.
  
  From (\ref{H5}) 
  \begin{eqnarray*}
    {\cal B}_{\infty} q &=& P_{\bmi{H}_2({\mathbb C}_0^+)}\left( G(\cdot) \hat{q}(-\cdot) \right) = P_{\bmi{H}_2({\mathbb C}_0^+)}\left( G(\cdot) G^{\dagger}(\cdot) \hat{v}(\cdot) \right) \\
    &=& P_{\bmi{H}_2({\mathbb C}_0^+)}\left( G(\cdot) \overline{G(\cdot)} \hat{v}(\cdot) \right) = P_{\bmi{H}_2({\mathbb C}_0^+)}\left(  \hat{v}(\cdot) \right)  =0,
    \end{eqnarray*}
    where we have used that $G$ is inner, and that $\hat{v} \in  \bmi{H}_2({\mathbb C}_0^-)$.
  
  For $u  \in \bmi{L}_2(0,\infty)$ with $u\perp \ker  {\mathcal B}_{\infty}$ we show next that 
  \begin{equation}
  \label{eq:16}
    \langle G(\cdot) \hat{u}(-\cdot), P_{\bmi{H}_2({\mathbb C}_0^-)}\left( G(\cdot) \hat{u}(-\cdot) \right) \rangle =0.
  \end{equation}
  We denote by $\hat{v}= P_{\bmi{H}_2({\mathbb C}_0^-)}\left( G(\cdot) \hat{u}(-\cdot)\right) $, and thus the inner product (\ref{eq:16}) becomes
  \begin{eqnarray*}
    \langle G(\cdot) \hat{u}(-\cdot), \hat{v}(\cdot) \rangle &=& \langle \hat{u}(-\cdot),\overline{G(\cdot)} \hat{v}(\cdot) \rangle\\
    &=& \langle \hat{u}(-\cdot),G^{\dagger}(\cdot) \hat{v}(\cdot) \rangle\\
    &=& \langle \hat{u}(\cdot),G^{\dagger}(-\cdot) \hat{v}(-\cdot) \rangle = \langle \hat{u}(\cdot), \hat{q}(\cdot) \rangle.
  \end{eqnarray*}
  The later is zero by the fact the inner product in Laplace domain equals the inner product in time domain, and that $u\perp \ker  {\mathcal B}_{\infty}$, $q \in \ker  {\mathcal B}_{\infty}$.
  
  Now we prove the equality (\ref{eq:15}). 
  \begin{eqnarray*}
    \| {\mathcal B}_{\infty} u \|^2 &=&  \langle P_{\bmi{H}_2({\mathbb C}_0^+)}\left( G(\cdot) \hat{u}(-\cdot) \right), P_{\bmi{H}_2({\mathbb C}_0^+)} \left(G(\cdot) \hat{u}(-\cdot) \right)  \rangle \\
   &=& \langle G(\cdot) \hat{u}(-\cdot) , P_{\bmi{H}_2({\mathbb C}_0^+)} \left(G(\cdot) \hat{u}(-\cdot) \right)  \rangle \\
   &=& \langle G(\cdot) \hat{u}(-\cdot) , G(\cdot) \hat{u}(-\cdot)  \rangle - \langle G(\cdot) \hat{u}(-\cdot) , P_{\bmi{H}_2({\mathbb C}_0^-)} \left(G(\cdot) \hat{u}(-\cdot) \right)  \rangle \\
     &=& \langle G(\cdot) \hat{u}(-\cdot) , G(\cdot) \hat{u}(-\cdot)  \rangle  - 0 =  \langle \hat{u}(-\cdot) , \hat{u}(-\cdot)  \rangle = \|u\|^2,
  \end{eqnarray*}
 where we have used (\ref{eq:16}) and isometry of Laplace and Fourier transforms, i.e. Paley-Wiener theorem.   
  From (\ref{eq:15}) it follows directly that the range of ${\mathcal B}_{\infty}$ is closed.
%
%
If we can show that the
  range is dense in $V$, then we have proved the assertion. Suppose
  that the range is not dense in $V$, then there exists a $v \in V$
  such that $v \in \left[{\rm Im}\, {\cal B}_{\infty}\right]^{\perp}$.
  So we have
  \begin{eqnarray*}
    0
    &=&
    \langle v, {\cal B}_{\infty}u \rangle\\
    &=&
    \frac{1}{2\pi} \int_{-\infty}^{\infty} v(i\omega)
    \overline{G(i\omega) \hat{u}(-i\omega)} d \omega\\
    &=&
    \frac{1}{2\pi} \int_{-\infty}^{\infty} \overline{G(i\omega)}
    v(i\omega) \,\overline{\hat{u}(-i\omega)} d \omega.
  \end{eqnarray*}
  Since this holds for every $\hat{u}(-\cdot) \in
  \bmi{H}_2({\mathbb C}_0^-)$, we have that
  \[
    x:=\overline{G} v \in \bmi{H}_2({\mathbb C}_0^+)
  \]
  Since $G$ is inner, we get from the above equation that
  \[
    v = G x,
  \]
  for the $x \in \bmi{H}_2({\mathbb C}_0^+)$. This means that $v
  \in V^{\perp}$, and since it is an element of $V$ it must be
  zero. This proves that the range of ${\cal B}_{\infty}$ equals
  $V$.
\item
  We see that it remains to prove that the semigroup has the
  additional properties.

  Gearhart \cite[Theorem 2.1 and Theorem 2.2]{ge78} or
  Moeller \cite[Theorem 3.1 and Theorem 3.2]{mo62} shows
  $\sigma(T(1))\subseteq \{ z\in\mathbb C\mid |z|<e^{-\alpha}\}$
  if and only if $\inf_{\re z\in (0,\alpha)} |G(z)| >0$.
  Thus this proves part 1.

  Moreover, from Gearhart \cite[Theorem 3.4]{ge78} we get that
  also part 2 of the theorem holds.
\hfill$\blacksquare$

\end{enumerate}

Next we give an example, which shows that Lemma \ref{lemH1} does
not hold without the assumption of infinite-time admissibility.
\begin{example}
Let $A_0$ be an infinitesimal generator on the infinite-dimensional
Hilbert space $H$ that satisfies $A_0=-A_0^*$, and let $b \in H$. Define
the operator $B_0 \in {\cal L}({\mathbb C},H)$ as $B_0u = b\cdot u$.
We assume that the system $(A_0, B_0^*)$ is approximately
observable.
It is well-known that
\[
  G(s):= 1 - B_0^*(sI - A_0 + \frac{1}{2} B_0B_0^*)^{-1}B_0
\]
is an inner function and Lemma \ref{lemH1} together with Theorem
\ref{theoreal} show that the Hankel operator $H_G$ has closed range.
Clearly, $(T_0, B_0, B_0^*,G)$ is a well-posed linear system,
where $T_0$ is the semigroup generated by $A_0-\frac{1}{2}B_0B_0^*$.
By the special structure of the system, we have that $(T_0, B_0,
B_0^*,G)$ is approximately controllable and approximately observable.
However, since $B_0$ is compact,
we have that $(A_0-\frac{1}{2}B_0B_0^*, B_0^*)$ is not
exactly observable in finite time.

We are now going to construct another realization.
We begin by considering the realization $(T,B,C,G)$ as constructed in
Theorem \ref{theoreal0} and Theorem \ref{theoreal}.
Note, that $(T,B,C,G)$ is exactly controllable,
exactly observable and the state space for the realization
in Theorem \ref{theoreal} is given by
\[
   V=[G \bmi{H}_2(\mathbb C_0^+)]^{\perp},
\]
where the orthogonal complement in taken in $\bmi{H}_2(\mathbb C_0^+)$.
We define $V_1$ as the closure of $V$ in the topology of
$\bmi{H}_2({\mathbb C}_{1}^+)$. Note that this space is isometric
isomorph (via Laplace transform) with the weighted $\bmi{L}_2$-space
\[
  \{ f \in \bmi{L}_2^{\rm loc}(0,\infty) \mid
  \int_0^{\infty} | e^{-t} f(t) |^2 dt < \infty \}.
\]
The semigroup on $V$ is given by
\[
  T(t) v = P_{\bmi{H}_2({\mathbb C}_0^+} ( e^{t\cdot} v ),
\]
It is easy to see that the inverse Laplace transform of $T(t)v$
is given by
\[
  \check{\left(T(t)v\right)}(\tau) = \check{v}(\tau + t), \qquad
  \tau \geq 0.
\]
Now we shall calculate the norm of $T(t)v$ in the new state space
$V_1$.
\begin{eqnarray*}
  \|T(t) v \|_{V_1}^2
  &=&
  \|\check{\left(T(t)v\right)} \|^2_{\check{V}_1}
  =
  \|\check{v}(\cdot + t) \|^2_{\check{V}_1}
  =
  \int_0^{\infty} |e^{-\tau} \check{v}(\tau + t) \|^2 d \tau\\
  &=&
  e^{2t} \int_0^{\infty} |e^{-\tau} \check{v}(\tau) \|^2 d \tau
  =
  e^{2t} \|v \|_{V_1}^2 .
\end{eqnarray*}
This implies that for every $t\geq 0$, we can extend the
operator $T(t)$ to $V_1$. Since $T(t)$ is a $C_0$-semigroup on $V$,
and since $V$ is dense in $V_1$, we have that the extension is
again a $C_0$-semigroup. We denote this new semigroup by $T_1(t)$.

Next we construct a well-posed realization of $G$ with state space
equal to $V_1$. As semigroup we take $T_1$, and as $B_1$ we take
$B_1 = G$.
 From part 2.\ in Theorem \ref{theoreal} we see that the
corresponding ${\cal B}_{1,\infty}$ satisfies
\begin{equation}
\label{admisB}
  {\cal B}_{1,\infty} = \imath {\cal B}_{\infty},
\end{equation}
where $\imath$ denotes the inclusion of $V$ into $V_1$. Thus
$B_1$ is an infinite-time admissible control operator for $T_1(t)$.
Furthermore, since the range of ${\cal B}_{\infty}$ equals $V$, and
since $V$ is dense in $V_1$, we have that $(T_1,B_1)$ is
approximately controllable.

As observation operator we take
\[
  C_1 x:=\check{x}(0),\quad x \in V_1.
\]
Hence the extension of $C$ to $V_1$.
We have to show that this is admissible observation operator for
$T_1(t)$. For $v \in V$, we have
that $C_1T_1(t) v = C T(t)v = \check{v}(t)$, and thus
\[
  \int_0^1 |\check{v}(t)|^2 dt
  \leq
  e^2 \int_0^1 |e^{-t} \check{v}(t) |^2 dt
  \leq
  e^2 \int_0^{\infty} |e^{-t} \check{v}(t) |^2 dt
  = e^2 \|v\|_{V_1}^2.
\]
Since $V$ is dense in $V_1$, we conclude that $C_1$ is an
admissible observation operator for $T_1(t)$. From
the definition of $C_1$ it is clear that $(T_1,C_1)$ is
approximately observable.

By (\ref{admisB}), we see that
\[
  \imath (zI-A)^{-1}B = (zI-A)^{-1} B_1,
\]
where $A_1$ is the infinitesimal generator of $T_1(t)$.
Hence
\begin{eqnarray*}
  C_1(sI-A_1)^{-1}(zI-A_1)^{-1}B_1
  &=&
  C_1(sI - A_1)^{-1}(zI-A)^{-1}B\\
  &=&
  C_1(sI - A)^{-1}(zI-A)^{-1}B\\
  &=&
  C(sI - A)^{-1}(zI-A)^{-1}B
  =\frac{G(s)-G(z)}{z-s}.
\end{eqnarray*}
Thus we have proved that $(T_1,B_1,C_1,G)$ is realization of $G$ as
well. Furthermore, this realization is approximately controllable
and observable.Thus the realizations $(T_0,B_0,B_0^*,G)$, $(T,B,C,G)$
and $(T_1,B_1,C_1,G)$
are all approximately controllable and approximately observable.

We now assume that all approximately controllable and
approximately observable realizations are equivalent.
The equivalence of $(T,B,C,G)$ and $(T_1,B_1$, $C_1,G)$ implies that the topologies
of $V$ and $V_1$ would be equivalent. Since $V_1$ is the closure
of $V$ in the topology of $V_1$, this implies that $V=V_1$.
In particular, there holds
\[
  \int_0^{\infty} e^{-2t} |f(t)|^2 dt \geq K \int_0^{\infty} |f(t)|^2
  dt,
\]
for all $f$ whose Laplace transform lies in $V_1$, where $K$ is independent
of $f$. Using this we see that
\begin{eqnarray*}
  K \int_0^{\infty} |f(t)|^2 dt
  &\leq&
  \int_0^{\infty} e^{-2t} |f(t)|^2 dt \\
  &=&
  \int_0^{t_0} e^{-2t} |f(t)|^2 dt + \int_{t_0}^{\infty} e^{-2t}
  |f(t)|^2 dt \\
  &\leq&
  \int_0^{t_0} |f(t)|^2 dt + e^{-2t_0} \int_{t_0}^{\infty} |f(t)|^2 dt \\
  &\leq&
  \int_0^{t_0} |f(t)|^2 dt + e^{-2t_0} \int_0^{\infty} |f(t)|^2 dt .
\end{eqnarray*}
Hence for all $t_0$
\[
  [K - e^{-2t_0} ] \int_0^{\infty} |f(t)|^2 dt \leq
  \int_0^{t_0} |f(t)|^2 dt
\]
Thus for $t_0$ sufficiently large
\begin{equation}
\label{ineq}
  \int_0^{t_0} |f(t)|^2 dt \geq K_1 \int_0^{\infty} |f(t)|^2 dt
\end{equation}
for all $f$'s whose Laplace transform lies in $V_1$.
For the output $y$ we have that
\[
  y(t) = C_1T_1(t) v = \check{v}(t)
\]
Hence with (\ref{ineq}), we obtain that
\[
  \int_0^{t_0} |y(t)|^2 dt \geq K_1 \int_0^{\infty} |y(t)|^2 dt
  =
  K_1 \int_0^{\infty} |\check{v}(t)|^2 dt = K_1 \|v\|^2_{V} \geq K_1
  \|v\|^2_{V_1}.
\]
Thus $(C_1,T_1)$ is exactly observable in finite-time.
The equivalence of the systems $(T_0,B_0,B_0^*,G)$ and
$(T_1,B_1,C_1,G)$
implies that the realization $(T_0,B_0,B_0^*,G)$
is also exactly observable in finite time.  However, this is not
possible, since $B_0^*$ is compact.

Concluding we see that the realization $(T_1,B_1,C_1,G)$ cannot be
equivalent with $(T,B,C,G)$ whereas $(T,B,C,G)$ is exactly controllable
and exactly observable, and $(T_1,B_1,C_1,G)$ is approximately controllable
and approximately observable. Note that the realization
$(T_1,B_1,C_1,G)$ does not have an infinite-time admissible
observation operator.
\end{example}

\section{Closing remarks}

The origin of this paper dates back to the research for the article \cite{JaZw02}. For that we needed that realization theory written down in the language of well-posed systems. Since that was not done before, we decided to do it ourselves. Hence we do not claim originality, but hope that this manuscript clarifies the ideas behind realization theory. For more reading we refer to Chapter 9 of \cite{Staf05}. Since there is a close link between properties of realizations and Hankel operators, the book of Peller \cite{Pell03} is also recommended.

\end{document}